# Taxicab Butterflies

by Kevin P. Thompson
Protective Life Insurance Company
Birmingham, AL  35223
thompson.kevin.p@gmail.com

In terms of the metric employed, taxicab geometry is a very close cousin to Euclidean geometry.  Instead of the $L^2$ metric, it utilizes the $L^1$ metric:

$$d\big((x_1, y_1), (x_2, y_2)\big) = |x_2 - x_1| + |y_2 - y_1|$$

This seemingly simple change in metric causes dramatic differences in the geometry: circles look like squares, π has the value 4, similar triangles all but vanish, and even distance from a point to a line becomes a more complicated concept [3, 4, 10].  The references provide a wide selection of books and articles that explore the nature of taxicab geometry and its relationship to Euclidean geometry.

This paper will examine the classic Euclidean geometry Butterfly Theorem in context of taxicab geometry.  While there are a number of proofs of the Euclidean Butterfly Theorem, many of them seem to leave the reader asking, "But, WHY is the theorem true?"  This investigation aims to provide a surprising level of insight into the core elements of the Euclidean version of the theorem.

**The Euclidean Butterfly Theorem**

The Butterfly Theorem is often formulated as follows: Let *M* be the midpoint of a chord *PQ* of a circle through which two other chords *AB* and *CD* pass, and let *X* be the intersection of *AB* and *PQ* and *Y* be the intersection of *CD* and *PQ*.  Then *M* is also the midpoint of *XY* (Figure 1).

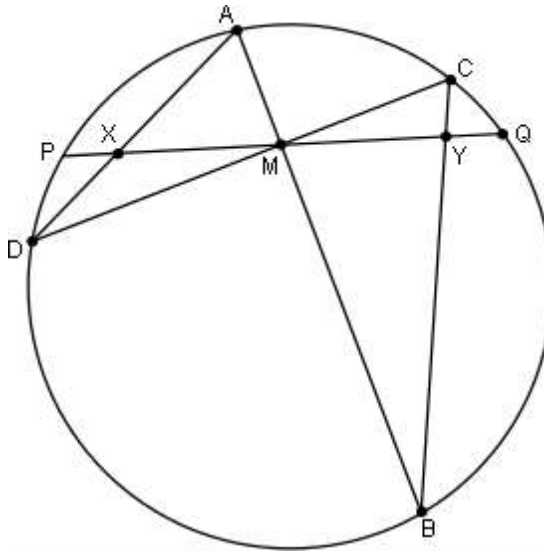

FIGURE 1: Classic illustration of the Euclidean Butterfly Theorem

The Euclidean Butterfly Theorem includes no restrictions.  For any choice of *P* and *Q* and any choice of *A* and *C*, the theorem is guaranteed to hold.

## The Taxicab Butterfly Theorem

First and foremost, the Butterfly Theorem simply does not hold in Taxicab geometry. Figure 1 illustrates an example where visually the theorem obviously fails to hold.

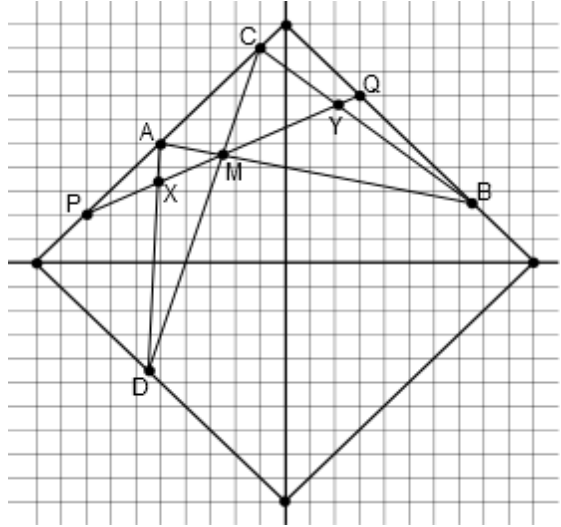

FIGURE 2: Failure of the Butterfly Theorem in Taxicab geometry

Quite a number of adjustments to the positions of $P$, $Q$, $A$, $C$, $B$, and $D$ can be made that fail to remedy the situation. In order to focus on the nature of the positive cases, the reader is invited to independently examine some of the additional negative cases summarized in the table below. All examples are based on a taxicab circle of radius 10.

| P | Q | A | C | Notes |
|---|---|---|---|---|
| (-8, 2) | (3, 7) | (-5, 5) | (1, 9) | $A$ and $C$ on the same segment of the circle |
| (-8, 2) | (3, 7) | (-2, 8) | (2, 8) | $A$ and $C$ symmetric about an axis of symmetry for the circle |
| (-5, 5) | (5, 5) | (-3, 7) | (2, 8) | $P$ and $Q$ symmetric about an axis of symmetry for the circle |
| (-9, 1) | (9, 1) | (-6, 4) | (-4, 6) | Points $A$ and $C$ and points $P$ and $Q$ symmetric about different axes of symmetry for the circle |
| $(-35/11, 75/11)$ | $(81/22, 139/22)$ | (-2, 8) | (2, 8) | Points $A$ and $C$ and points $B$ and $D$ symmetric about different axes of symmetry for the circle |

TABLE 1: Additional examples where the Butterfly Theorem fails to hold in taxicab geometry

From the examples above, it appears that a high degree of symmetry is necessary for the Butterfly Theorem to hold in taxicab geometry. And, this is indeed the case. Figures 3, 4, and 5 illustrate situations in which the theorem holds. In all of the successful examples, the points $A$ and $C$ and the points $B$ and $D$ are symmetric on the taxicab circle in the same manner, meaning both pairs are symmetric about an axis of symmetry of the circle or both pairs are symmetric about the center of the circle.

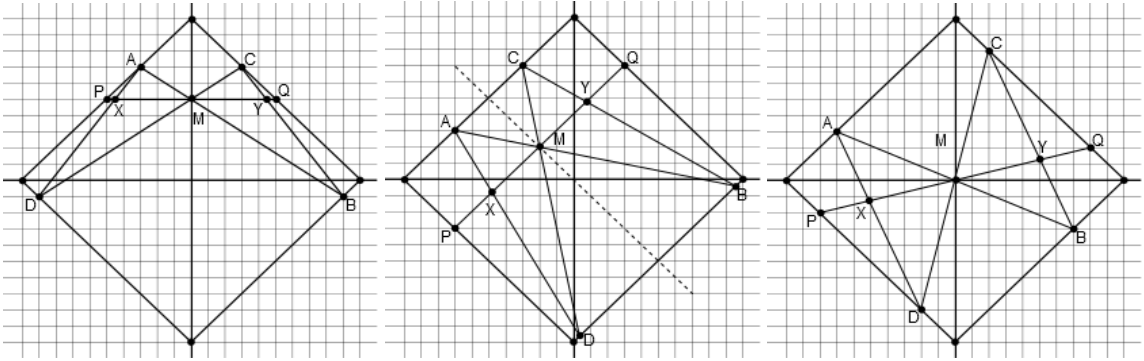

FIGURES 3–5: Successful demonstrations of the Butterfly Theorem in Taxicab geometry

The examples provided thus far offer enough insight to formulate the taxicab geometry version of the Butterfly Theorem which can be proved with simple symmetry arguments.

THEOREM (Taxicab Butterfly Theorem): Let *M* be the midpoint of a chord *PQ* of a circle through which two other chords *AB* and *CD* pass, and let *X* be the intersection of *AB* and *PQ* and *Y* be the intersection of *CD* and *PQ*. Then *M* is also the midpoint of *XY* when the points *A* and *C* and the points *B* and *D* are symmetric on the circle in the same manner – either about an axis of symmetry or about the center of the circle.

When *A* and *C* and *B* and *D* are symmetric on the circle, conditions are implicitly placed on the position of the midpoint *M* of *PQ* (and *XY*). Therefore, an alternate (but equivalent) formulation of the theorem might in some cases be easier to understand or apply.

THEOREM (Taxicab Butterfly Theorem, Alternate): Let *M* be the midpoint of a chord *PQ* of a circle through which two other chords *AB* and *CD* pass, and let *X* be the intersection of *AB* and *PQ* and *Y* be the intersection of *CD* and *PQ* (Figure 1). Then *M* is also the midpoint of *XY* when either:
  1. *M* is the center of the circle, OR
  2. *A* and *C* are symmetric on the circle and *M* lies on the axis of symmetry for *A* and *C*

Note that the conditions on *A*, *B*, *C*, and *D* in the primary version and the conditions on *M* in the alternate version also force *P* and *Q* to be symmetric on the taxicab circle.

## Shedding Light on the Euclidean Butterfly Theorem

The first formulation of the Taxicab Butterfly Theorem is probably the most useful in giving insight to the Euclidean Butterfly Theorem. In this formulation, the geometries are in perfect agreement – the points *A* and *C* and the points *B* and *D* are symmetric on the circle in the same manner. In Euclidean geometry, there is always an axis of symmetry between any two points on the circle. Despite the most common examples of proof, this appears to be the key ingredient to the theorem. Since Euclidean circles are "abundantly symmetric" (every possible diameter is an axis of symmetry) the Euclidean theorem holds for any choice of *P*, *Q*, *A*, and *C*. And, this is the downfall for the theorem in taxicab geometry – taxicab circles have precisely 4 axes of symmetry.

For the second formulation, there is partial agreement between the geometries. While *A* and *C* are always symmetric on a Euclidean circle, there is no requirement for *M* to either be the center of the circle or to lie on the axis of symmetry for *A* and *C* (see Figure 1 or Figure 6). In taxicab geometry, this condition is a consequence of the symmetry requirements for *A*, *C*, *B*, and

*D* and is likely necessary to gain sufficient symmetry from the taxicab circle to force a symmetry relationship between *A–C* and *M*. Euclidean circles have so much symmetry that there appears to be an implied symmetry relationship between *A–C* and *M* that must be made explicit in the taxicab geometry case due to its relative lack of symmetry.

Other aspects to consider are the concepts of similar and congruent triangles. Similar triangles are the crux of many of the proofs of the Euclidean Butterfly Theorem. But, similar triangles are much harder to come by in taxicab geometry. Two triangles must satisfy a SASAS condition in order to conclude similarity [10].

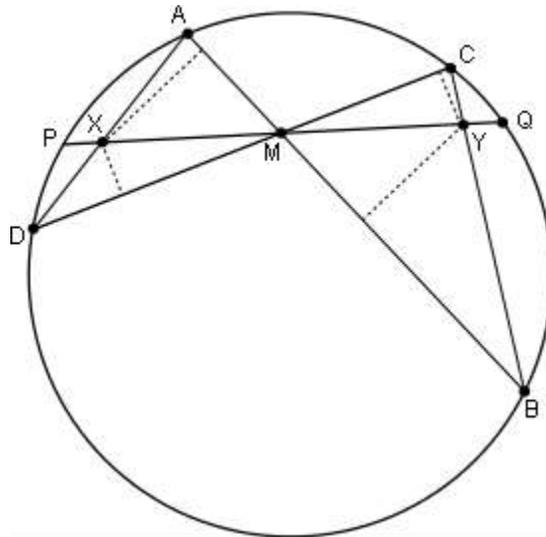

FIGURE 6: Detailed illustration of the Butterfly Theorem in Euclidean geometry

Very much related to the struggle for similar triangles in taxicab geometry is the lack of preservation of length under rotation for the taxicab metric. Consider the four central triangles in Figure 6. One of the conclusions of the Euclidean Butterfly Theorem is that two pairs of these triangles are congruent and thus are reflections about the point *M* (using *PQ* and its perpendicular as axes for reflection). We could rotate this circle by any degree and the result will remain true. But, this is not the case in taxicab geometry. Only 4 precise rotations will preserve the taxicab length of a line (2, 4, 6, and 8 t-radians), and only reflections about lines parallel to the *x*-axis, *y*-axis, and the lines $y = x$ and $y = -x$ will preserve the taxicab length of a line [1].

If indeed one of the basic conditions in the Butterfly Theorem is for certain triangles to be symmetric about the point *M*, then the situations in which the theorem holds in taxicab geometry will be limited to those instances where the axes of reflection about *M* are parallel to the 4 axes of symmetry of the taxicab circle, which is what the examples and discussion above seem to imply.

**Conclusion**

The Taxicab Butterfly Theorem appears to hold only in the narrow circumstances where the points *A* and *C* and the points *B* and *D* are symmetric on the taxicab circle in the same manner – either about an axis of symmetry of the circle or about the center of the circle.

While the true nature and meaning of the Euclidean Butterfly Theorem still seems somewhat elusive, it does appear that abundant symmetry in the Euclidean circle plays a crucial role. In addition, the relatively lenient requirements for Euclidean similar triangles (when compared to taxicab geometry) and the preservation of length under rotation in the Euclidean metric also provide some insight as to why the theorem holds under any conditions in Euclidean geometry.

To learn more about taxicab geometry, visit the Taxicab Geometry website at http://www.taxicabgeometry.net